\documentclass[11pt]{article}
\usepackage{amsmath,amsfonts,amsthm,amssymb}
\usepackage{url}
\usepackage{graphicx}
\usepackage{lmodern}
\title{A note on partial rejection sampling for the hard disks model in the plane}
\author{Jake Wellens}

\usepackage{color}

\usepackage[colorlinks=true, linkcolor=red, urlcolor=blue, citecolor=blue]{hyperref}

\graphicspath{ {images/} }

\newcommand\be{\begin{equation}}
\newcommand\ee{\end{equation}}
\newcommand\bea{\begin{eqnarray}}
\newcommand\eea{\end{eqnarray}}

\newcommand\R{\mathbb{R}}

\newcommand\E{\mathbb{E}}

\newcommand\nn{\nonumber}

\newcommand{\BadPairs}{\text{\fontfamily{lmss}\selectfont
BadPairs}}
\newcommand{\BadPoints}{\text{\fontfamily{lmss}\selectfont
BadPoints}}
\newcommand{\Area}{\text{Area}}

\newtheorem{thm}{Theorem}

\newtheorem{lem}[thm]{Lemma}

\newtheorem{fact}[thm]{Fact}

\advance \topmargin by -\headheight
\advance \topmargin by -\headsep
\evensidemargin \oddsidemargin
\begin{document}
\maketitle
\begin{abstract}

In this note, we slightly improve the guarantees obtained by Guo and Jerrum for sampling from the hard disks model in the plane via partial rejection sampling. Our proof makes use of the fact that if one spreads apart a collection of disks in the plane, the area of the union of the disks cannot decrease.

\end{abstract}

\section{Introduction}

Coming from statistical physics, the \emph{hard disks model} is a simple probability distribution used to model positions of particles of a contained gas, supported on configurations non-overlapping disks of radius $r$ in a bounded region of $\R^2$. More precisely, the centers of the disks are sampled from a Poisson process of intensity $\lambda_r = \lambda/(\pi r^2)$, conditioned on the disks being non-overlapping. 

In a very recent work of Guo and Jerrum \cite{GJ}, the general-purpose method of \emph{partial rejection sampling} is applied to the problem of sampling from the hard disks model. The authors show that for sufficiently small $\lambda > 0$, this algorithm is efficient -- namely it runs in expected $O(r^{-2})$ time. This bound essentially follows from the following theorem:

\begin{thm}\label{GJ hard disks}\emph{(Guo, Jerrum \cite{GJ})} Partial rejection sampling for the hard disks model with parameter $\lambda_r = \frac{\lambda}{\pi r^2}$ takes $O(\log{1/r})$ rounds of resampling in expectation for $\lambda \leq \overline{\lambda} = 0.21027$. 
\end{thm}

The authors of \cite{GJ} conjecture that $\overline{\lambda}$ can be taken to be $\approx 0.5$, according to their simulations. In this note, we slightly improve the value of $\overline{\lambda}$ (although it remains quite far from the conjectured truth):

\begin{thm}\label{improved}The constant $\overline{\lambda}$ in Theorem \ref{GJ hard disks} can be improved to $0.2344+$.
\end{thm}

Our analysis closely resembles that given in \cite{GJ}, with an extra ingredient. A theorem of Bollobas \cite{Bollobas} states that if one pushes apart a collection of disks in $\R^2$ (of equal radii) in a continuous fashion such that the pairwise distances between their centers are always increasing, then the area of their union is also increasing. In particular, we have the following special case:

\begin{fact}\label{bollobas} Let $\gamma > 1$, and $x_1, \dots, x_n \in \R^2$. Then
$$\emph{\Area}\left(\bigcup_{i=1}^n D_1(x_i)\right) \leq \emph{\Area}\left(\bigcup_{i=1}^n D_1(\gamma x_i)\right).$$
\end{fact}

While this may seem intuitively obvious, proving the statement for general expansions of disks with different radii was an open problem until 2002, and it is still open in dimensions higher than 2. In any case, with this fact in hand, we now prove Theorem \ref{improved}.

\section{Proof of Theorem \ref{improved}} 

The improvement over the estimate in \cite{GJ} mostly boils down to the following consequence of Fact \ref{bollobas}:

\begin{lem}\label{geom}
Let $C = \cup_{i=1}^{\ell} D_{2r}(x_i)$ be the union of $\ell$ disks of radius $2r$ in $\R^2$. Then 
$$\frac{1}{2}\int_C \int_C 1_{\|x - y\| \leq 2r} \, dx \, dy \geq \frac{1}{2} \frac{\emph{\Area}(C)}{4\pi r^2} \cdot \int_{D_{2r}(0)} \int_{D_{2r}(0)} 1_{\|x - y\| \leq 2r} \, dx \, dy $$
\end{lem}

\textit{Proof:} Set $f(x) := \Area(C \cap D_{2r}(x))$ and $g(x) := \Area(D_{2r}(0) \cap D_{2r}(x))$. Then proving the lemma is equivalent to showing that
\be
\frac{1}{\Area(C)} \int_C f(x) \, dx \geq \frac{1}{4\pi r^2} \int_{D_{2r}(0)} g(x) \, dx
\ee
which is in turn equivalent to the inequality
\be
\frac{1}{\Area(C)} \int_{0}^{\infty} \Area(x \in C: f(x) > t)  \, dt \geq \frac{1}{4\pi r^2} \int_{0}^{\infty} \Area(x \in D_{2r}(0): g(x) > t) 
\ee
We will show something even stronger: for each $t \geq 0$, we have
\be \label{stronger}
\frac{ \Area(x \in C: f(x) > t)}{\Area(C)} \geq \frac{\Area(x \in D_{2r}(0): g(x) > t)}{4\pi r^2}
\ee
Observe that for each $t$, the set $\{x \in D_{2r}(0): g(x) > t\}$ is an open disk $D_{2\alpha r}(0)$ for some $\alpha = \alpha(t) \leq 1$. Then clearly $\cup_{i=1}^{\ell} D_{2\alpha r}(x_i) \subseteq \{x \in C: f(x) > t\}$. Hence to prove (\ref{stronger}), it suffices to show that
\be \label{pre-scale}
\frac{\Area(\cup_{i=1}^{\ell} D_{2\alpha r}(x_i))}{\pi(2\alpha r)^2} \geq \frac{\Area(\cup_{i=1}^{\ell} D_{2r}(x_i)) }{\pi(2r)^2}.
\ee
Consider applying the transformation $x \mapsto x/2\alpha r$ on $\R^2$. This sends $D_{2\alpha r}(x_i) \mapsto D_1(x_i/2\alpha r)$ and scales all areas by $1/(2\alpha r)^2$. Similarly the transformation $x \mapsto x/2r$ takes $D_{2r}(x_i) \mapsto D_{1}(x_i/2r)$ and scales areas by $1/(2r)^2$. Hence, (\ref{pre-scale}) is equivalent to
\be
\Area\left(\bigcup_{i=1}^\ell D_1(x_i/2\alpha r)\right) \geq \Area\left(\bigcup_{i=1}^\ell D_1(x_i/2r)\right)
\ee
which is Fact \ref{bollobas}, with $\gamma = 1/\alpha$. \qed

\textit{Proof of Theorem \ref{improved}:} We use the same notation as \cite{GJ}. Recall that $P_t$ is the set of points sampled during round $t \geq 0$ of the algorithm. Suppose $|\BadPairs(P_t)| = k_t$, where $\BadPairs(P_t)$ is the set of \emph{unordered} pairs $\{x, y\} \subset P_t$ with $\|x-y\| < 2r$, and $\BadPoints(P_t)$ is the set of all points which occur in a bad pair. Then the resampling set $S_t \subset[0,1]^2$ is $\cup_{x \in \BadPoints(P_t)} B_{2r}(x)$. Let
\bea \nn
k' &=& \E[k_{t+1}\, | \, \BadPoints(P_t)] \\ \nn
j' &=& \E[\#\{(x,y) \in S_t \times [0,1]^2: \|x-y\| \leq 2r\}\cap P_{t+1}^2 \, | \, \BadPoints(P_t)] \\ \nn
\ell' &=& \E[\#\{(x,y) \in S_t \times S_t : \|x-y\| \leq 2r\} \cap P_{t+1}^2\, | \, \BadPoints(P_t)]
\eea
Since each unordered bad pair $\{x,y\}$ with $x, y \in S_t$ gets counted twice in $j'$ and twice in $\ell'$, we have
 \be
 k' = j' - \frac{\ell'}{2}.
 \ee
It is shown in \cite{GJ} that the hard disks process can be coupled to a Poisson process in such a way that the latter configuration always contains the former. In particular,
\bea
j' &\leq& \int_{S_t} \lambda_r \int_{[0,1]^2}  \lambda_r 1_{\|x-y\|\leq 2r} \, dx\, dy \\ &\leq& \lambda_r^2 \cdot \Area(S_t) \cdot 4 \pi r^2 = \frac{4\lambda^2}{\pi r^2} \cdot \Area(S_t)
\eea
On the other hand, by Lemma \ref{geom}, we have
\bea
\frac{\ell}{2} &=& \frac{1}{2} \int_{S_t} \lambda_r \int_{S_t}  \lambda_r 1_{\|x-y\|\leq 2r} \, dx\, dy \\
 &\geq& \frac{\Area(S_t)}{4\pi r^2} \cdot \underbrace{\frac{1}{2} \int_{D_{2r}(0)} \lambda_r \int_{D_{2r}(0)}  \lambda_r 1_{\|x-y\|\leq 2r} \, dx\, dy}_{= (8 - \frac{6\sqrt{3}}{\pi})\lambda^2}
\eea
(see \cite{GJ} for an evaluation of the integral) and hence
\be
k' \leq \Area(S_t) \cdot \lambda^2 \cdot\underbrace{\left(\frac{4}{\pi r^2} - \frac{2}{\pi r^2} + \frac{3\sqrt{3}}{2\pi^2r^2}\right)}_{=\frac{4\pi+3\sqrt{3}}{2\pi^2r^2}}
\ee
Recall that $S_t$ is the union of disks of radius $2r$ centered at $\leq 2k_t$ points, each of which must overlap with at least one other disk -- in fact, for each disk $D_{2r}(x)$, there must be an overlapping disk $D_{2r}(x')$ with $\|x - x'\| < 2r$. Therefore $\Area(S_t)$ is maximized when $S_t$ is a union of $k_t$ connected components, each of which is a translated copy of $D_{2r}(0,0)\cup D_{2r}(0,2r)$, which has area $(\frac{16\pi}{3} + 2\sqrt{3})r^2$, as can be seen by elementary geometry. Finally we obtain the estimate 
\be
k' \leq \left(\frac{16\pi}{3} + 2\sqrt{3}\right)\cdot \frac{4\pi+3\sqrt{3}}{2\pi^2} \cdot \lambda^2 \cdot k_t 
\ee
and so we may take $\overline{\lambda} = 0.2344+$. \qed \\

\end{document}